\documentclass[12pt]{article}
\usepackage{amsmath}
\usepackage{amssymb}
\usepackage{amscd}
\usepackage{amsthm}

\theoremstyle{plain}
\newtheorem{Thm}{Theorem}

\newtheorem{Lem}[Thm]{Lemma}

\theoremstyle{definition}

\numberwithin{equation}{section}

\title{Flops connect minimal models}
\author{Yujiro Kawamata}

\begin{document}

\maketitle

\begin{abstract}
A result by Birkar-Cascini-Hacon-McKernan together with the
boundedness of length of extremal rays implies that different
minimal models can be connected by a sequence of flops.
\end{abstract}

A {\em flop} of a pair $(X,B)$ is a flip of a pair $(X,B')$ which is
crepant for $K_X+B$ where $B'$ is a suitably chosen different
boundary. We prove the following:

\begin{Thm}
Let $f: (X,B) \to S$ and $f': (X',B') \to S$ be projective morphisms
from $\mathbf{Q}$-factorial terminal pairs of varieties and
$\mathbf{Q}$-divisors
such that $K_X+B$ and $K_{X'}+B'$ are relatively
nef over $S$.
Assume that there exists a birational map $\alpha: X \dashrightarrow X'$
such that $\alpha_*B=B'$, where the lower asterisk denotes the strict transform.
Then $\alpha$ is decomposed into a sequence of flops.
\end{Thm}

More precisely, there exist an effective $\mathbf{Q}$-divisor $D$ on
$X$ such that $(X,B+D)$ is klt and a factorization of the birational
map $\alpha$
\[
X = X_0 \dashrightarrow X_1 \dashrightarrow \cdots \dashrightarrow X_t
= X'
\]
which satisfy the following conditions:

(1) $\alpha_i: X_{i-1} \to X_i$ ($1 \le i \le t$) is a flip for the pair
$(X_i,B_i+D_i)$ over $S$, where $B_i$ and $D_i$ are strict
transforms of $B$ and $D$, respectively.

(2) $\alpha_i$ is crepant for $K_{X_{i-1}}+B_{i-1}$ in the sense that the
pull-backs of $K_{X_{i-1}}+B_{i-1}$ and $K_{X_i}+B_i$ coincide on a
common log resolution.

We remark that the boundary $B$ need not be assumed to be big as in
\cite{BCHM}~Corollary~1.1.3. For example, a birational map between
Calabi-Yau manifolds can be decomposed into a sequence of flops. The
number of {\em marked} minimal models which are birationally equivalent to
a fixed pair is
finite if $B$ is big (\cite{BCHM}~Corollary~1.1.5), but it is not
the case in general (cf. \cite{CY}), where a marked minimal model is
a pair consisting of a minimal model and a fixed birational map to
it. If we relax the condition for the pairs to being klt, then we
should allow crepant blowings up besides flops.

The theorem was already proved in the case $\dim X = 3$ and $B=0$;
first in \cite{crepant} assuming the abundance which was proved afterwards,
and later in \cite{Kollar} without assumption.

\begin{proof}
It is well-known that $\alpha$ is an isomorphism in codimension $1$
because $(X,B)$ and $(X',B')$ are terminal and $K_X+B$ and
$K_{X'}+B'$ are relatively nef (cf. \cite{crepant}).
We recall the proof for reader's convenience.
Let $\mu: V \to X$ and $\mu': V \to X'$ be common log resolutions.
We write
\[
K_V = \mu^*(K_X+B) - \mu_*^{-1}B + E = (\mu')^*(K_{X'}+B') -
(\mu'_*)^{-1}B' + E'
\]
where $E$ and $E'$ are effective divisors whose supports coincide
with the exceptional loci of $\mu$ and $\mu'$, respectively, 
because $(X,B)$ and $(X',B')$ are
terminal. Assume that there is a prime divisor on $V$ which is
contracted by $\mu$ but not by $\mu'$.  Then it is an irreducible
component of $E$ but not of $E'$. We set $F = \min \{E,E'\}$, $\bar
E = E - F$ and $\bar E' = E' - F$. By the Hodge index theorem, there
exists a curve $C$ on $V$ which is contracted by $\mu$ and is
contained in $\text{Supp}(\bar E)$ but not in
$\text{Supp}(\mu_*^{-1}B+\bar E')$ and such that $(\bar E \cdot C)
<0$. Since $\mu_*^{-1}B \ge (\mu'_*)^{-1}B'$, we have
\[
((\mu')^*(K_{X'}+B') + \mu_*^{-1}B - (\mu'_*)^{-1}B' + \bar E')
\cdot C) \ge 0.
\]
But this is a contradiction to
\[
(\mu^*(K_X+B) + \bar E) \cdot C) < 0.
\]
The case where there is a prime divisor on $V$ which is
contracted by $\mu'$ but not by $\mu$ is treated similarly.

Let $L'$ be an effective $f'$-ample divisor on $X'$,
and $L$ its strict transform on $X$.
There exists a small positive number $l$ such that
$(X,B+lL)$ is klt.
If $K_X+B+lL$ is $f$-nef over $S$, then
$\alpha$ becomes a morphism by the base point free theorem,
hence an isomorphism since $X$ is $\mathbf{Q}$-factorial.
Therefore we may assume that $K_X+B+l'L$ is not
$f$-nef over $S$ for any $0 < l' \le l$.

Let $H$ be an effective divisor on $X$ such that $(X,B+lL+tH)$ is
klt and $K_X+B+lL+tH$ is $f$-nef for some positive number $t$. We
shall run the MMP for the pair $(X,B+l'L)$ over $S$ with scaling of
$H$ for some $l'$. Since $\alpha$ is an isomorphism in codimension
$1$, there are only flips in this MMP. The following lemma shows
that we can choose extremal rays such that the flips are crepant
with respect to $K_X+B$.

Let $k$ be a positive integer such that
$k(K_X+B)$ is a Cartier divisor.
We set $e = \frac 1{2k\dim X+1}$.

\begin{Lem}
(1) There exists an extremal ray $R$ for $(X,B+lL)$ over $S$ such
that $((K_X+B) \cdot R)=0$.

(2) Let
\[
\begin{split}
t_0 = \min &\{t \in \mathbf{R} \,\vert\,
((K_X+B+lL +tH) \cdot R) \ge 0 \text{ for all extremal rays } R \\
&\text{ for } (X, B + lL) \text{ over } S \text{ s.t. }
((K_X+B) \cdot R ) = 0\}.
\end{split}
\]
Then $K_X+B+elL + et_0H$ is $f$-nef, and there exists an extremal
ray $R$ for $(X, B + elL)$ over $S$ such that $((K_X+B+elL + et_0H)
\cdot R) = ((K_X+B) \cdot R)=0$.
\end{Lem}

\begin{proof}
(1) Since $K_X+B+elL$ is not nef, there exists an extrenal ray $R$ for
$(X, B+elL)$ over $S$.
Then $R$ is also an extremal ray for $(X,B+lL)$ because $(X,B)$ is $f$-nef.
Since the pair $(X,B+lL)$ is klt, 
$R$ is generated by a rational curve $C$, 
which is mapped to a point on $S$, such that
\[
0 > ((K_X+B+lL) \cdot C) \ge -2\dim X
\]
by \cite{length}.

We claim that $((K_X+B) \cdot C)=0$.
Indeed we have otherwise $((K_X+B) \cdot C) \ge 1/k$, hence
\[
\begin{split}
&((K_X+B+elL) \cdot C) \\
&= \frac 1{2k\dim X+1} ((K_X+B+lL) \cdot C)
+ \frac{2k\dim X}{2k\dim X+1}((K_X+B) \cdot C) \\
&\ge \frac 1{2k\dim X+1}(-2\dim X+2\dim X) = 0
\end{split}
\]
a contradiction. 

(2) If $K_X+B+elL + et_0H$ is not $f$-nef, then there exists an
extremal ray $R$ for $(X,B+elL + et_0H)$ over $S$. Then $R$ is also
an extremal ray for $(X,B+lL+t_0H)$ because $(X,B)$ is $f$-nef.
Since the pair $(X,B+lL+t_0H)$ is klt, $R$ is generated by a
rational curve $C$ such that $((K_X+B+lL+t_0H) \cdot C) \ge -2\dim
X$ by \cite{length}. Then we have
\[
\begin{split}
&((K_X+B+elL+et_0H) \cdot C) \\
&= \frac 1{2k\dim X+1} ((K_X+B+lL+t_0H) \cdot C)
+ \frac{2k\dim X}{2k\dim X+1}((K_X+B) \cdot C) \\
&\ge \frac 1{2k\dim X+1}(-2\dim X+2\dim X) = 0
\end{split}
\]
a contradiction.  Therefore $K_X+B+elL + et_0H$ is $f$-nef.

Since $B+lL$ is $f$-big, the number of extremal rays for $(X,B+lL)$ over $S$ 
is finite.
Hence there exists such an $R$ that 
$((K_X+B+lL +t_0H) \cdot R)=((K_X+B) \cdot R) = 0$.
\end{proof}

We note that we can deduce (1) from only the finiteness of extremal
rays, but not (2).
The point is that the number $e$ stays independent of $t_0$ during the MMP.

We run the MMP for $(X,B+elL)$ with scaling of $H$. We take an
extremal ray $R$ such that $((K_X+B+elL + et_0H) \cdot R) = ((K_X+B)
\cdot R)=0$. The flip exists by \cite{BCHM}~Corollary~1.4.1. Since
$(lL + t_0H) \cdot R)=0$, the pair $(X,B+lL+t_0H)$ remains to be klt
after the flip. We also note that $k(K_X+B)$ remains to be a Cartier
divisor after the flip by the base point free theorem.
Therefore we can continue the process. By the termination theorem of
directed flips (\cite{BCHM}~Corollary~1.4.2), we complete our proof.
\end{proof}

Department of Mathematical Sciences, University of Tokyo,

Komaba, Meguro, Tokyo, 153-8914, Japan

kawamata@ms.u-tokyo.ac.jp


\begin{thebibliography}{length}

\bibitem{BCHM}
Caucher Birkar, Paolo Cascini, Christopher D. Hacon, James McKernan.
{\em Existence of minimal models for varieties of log general type}.
math.AG/0610203.

\bibitem{crepant}
Kawamata, Yujiro.
{\em Crepant blowing-up of 3-dimensional canonical singularities
and its application to degenerations of surfaces}.
Ann. of Math. {\bf 127} (1988), 93--163.

\bibitem{length}
Kawamata, Yujiro.
{\em On the length of an extremal rational curve}.
Invent. Math. {\bf 105} (1991), 609--611.

\bibitem{CY}
Kawamata, Yujiro.
{\em On the cone of divisors of Calabi-Yau fiber spaces}.
Internat. J. Math. {\bf 8} (1997), 665--687.

\bibitem{Kollar}
Koll\'ar, J\'anos.
{\em Flops}.
Nagoya Math. J. {\bf 113}(1989), 15--36.


\end{thebibliography}
\end{document}